\documentclass[11pt]{article}

\usepackage[cp850]{inputenc}
\usepackage{latexsym,graphicx}

\font\tenmath=msbm10 scaled 1200
\font\sevenmath=msbm7 scaled 1200
\font\fivemath=msbm5 scaled 1200

\newfam\mathfam \textfont\mathfam=\tenmath
\scriptfont\mathfam=\sevenmath \scriptscriptfont\mathfam=\fivemath
\def\math{\fam\mathfam}
\def\R{{\math R}}
\def\N{{\math N}}
\def\E{{\math E}}

\def\P{{\math P}}

\def \^#1{\if#1i{\accent"5E\i}\else{\accent"5E#1}\fi}

\def \lesim{\stackrel{<}{\sim}}

\def \g{\gamma}

\def \ni{\noindent}

\newtheorem{Thm}{Theorem}

\newtheorem{Lem}{Lemma}
\newtheorem{Pro}{Proposition}
\newtheorem{Cor}{Corollary}

\author{\sc
{\sc Harald Luschgy}\thanks{Universit\"at Trier, FB IV-Mathematik, D-54286 Trier, BR Deutschland.
E-mail: {\tt luschgy@uni-trier.de}}
\quad {and}
\quad {\sc Gilles Pag\`es} \thanks{Laboratoire de Probabilit\'es et Mod\`eles al\'eatoires, UMR~7599, Universit\'e Paris 6, case 188, 4,
pl. Jussieu, F-75252 Paris Cedex 5. E-mail:{\tt  gpa@ccr.jussieu.fr}}
}

\title{\bf  Functional quantization and metric entropy for Riemann-Liouville processes}

\begin{document}

\maketitle
\begin{abstract}
We derive a high-resolution formula for the $L^2$-quantization errors of Riemann-Liouville processes and the sharp Kolmogorov 
entropy asymptotics for related Sobolev balls.
We describe a quantization procedure which leads to asymptotically optimal functional quantizers. Regular variation of 
the eigenvalues of the covariance operator plays a crucial role.
\end{abstract}


\bigskip
\ni {\em MSC: 60G15, 60E99, 41A46.}

\section{Introduction}
\setcounter{equation}{0}
\setcounter{Assumption}{0}
\setcounter{Theorem}{0}
\setcounter{Proposition}{0}
\setcounter{Corollary}{0}
\setcounter{Lemma}{0}
\setcounter{Definition}{0}
\setcounter{Remark}{0}

Functional quantization of stochastic processes can be seen as a discretization of the path-space of a process and the
 approximation (coding) of a process by finitely many deterministic functions from its path-space. In a Hilbert space setting this reads as follows.

Let $(H, < \cdot , \cdot >) $ be a separable Hilbert space with norm
$\|  \cdot \| $
and let $X : (\Omega, {\cal A}, \P) \rightarrow H$ be a random vector taking its values in $H$ with distribution $\P_X$. For $n \in \N$, the
$L^2$-quantization problem for $X$ of level $n$ (or of nat-level $\log n$) consists in minimizing
\[
\left(\E \min_{a \in \alpha} \| X - a \|^2 \right)^{1/2} =
\|
\min_{a \in \alpha} \| X-a \|  \|_{L^2(\P)}
\]
over all subsets $\alpha \subset H$ with $\mbox{card} (\alpha) \leq n$. Such a set $\alpha$ is called $n$-codebook or $n$-quantizer. 
The minimal $n$th quantization
error of $X$ is then defined by
\begin{equation}
e_n (X) := \inf \left\{ ( \E \min_{a \in \alpha} \| X-a \|^2)^{1/2} : 
\alpha \subset H, \; \mbox{card} (\alpha)
\leq n \right\} .
\end{equation}
Under the integrability condition
\begin{equation}\label{(1.2)}
\E \,\|  X \|^2 < \infty
\end{equation}
the quantity $e_n (X)$ is finite.

For a given $n$-quantizer $\alpha$ one defines an associated closest neighbour projection
\[
\pi_\alpha := \sum\limits_{a \in \alpha} a \mbox{\bf 1}_{C_{a} (\alpha)}
\]
and the induced $\alpha$-quantization (Voronoi quantization) of $X$ by
\begin{equation}
\hat{X}^\alpha := \pi_\alpha (X) ,
\end{equation}
where $\{ C_a (\alpha) : a \in \alpha\}$ is a Voronoi partition induced by $\alpha$, that is a Borel partition of $H$ satisfying
\begin{equation}
C_a (\alpha) \subset V_a (\alpha) := \{ x \in H : \| x-a \| = \min_{b \in \alpha} \| x-b \|  \}
\end{equation}
for every $a \in \alpha$. Then one easily checks that, for any random vector $X^{'} : \Omega \rightarrow \alpha \subset H$,
\[
\E\,\| X - X^{'} \|^2 \geq \,\E\,\| X - \hat {X}^\alpha \|^2 = \E\,\min_{a \in \alpha} \| X-a \|^2
\]
so that finally
\begin{eqnarray}
e_n (X) & = & \inf \left\{ ( \E\,\| X - \hat{X} \|^2)^{1/2} : \hat{X} = f(X), f : H \rightarrow H \;
\mbox{Borel measurable,} \; \right.\\
             && \left.\qquad \mbox{card} (f(H)) \leq n \right\} \nonumber \\
        & = & \inf \left\{ ( \E\,\| X- \hat{X} \|^2)^{1/2} : \hat{X} : \Omega \rightarrow H \;
\mbox{random vector,} \;
             \mbox{card} (\hat{X}(\Omega)) \leq n \right\} . \nonumber
\end{eqnarray}
Observe that the Voronoi cells $V_a (\alpha), a \in \alpha$ are closed and convex (where convexity is a characteristic 
feature of the underlying Hilbert structure). 
Note further that there are infinitely many $\alpha$-quantizations of $X$ which all produce the same quantization 
error and $\hat{X}^\alpha$
is $\P$-a.s. uniquely defined if $\P_X$ vanishes on hyperplanes.

A typical setting for functional quantization is $H = L^2([0,1], dt)$ but is obviously not restricted to the Hilbert space setting.
Functional quantization is the natural extension to stochastic processes of the so-called optimal vector
quantization of random vectors in $H = \R^d$ which has been extensively investigated since the late 1940's in Signal 
processing and Information Theory (see \cite{GERSH}, \cite{GRAY}).
For the mathematical aspects of vector quantization in $\R^d$, one may consult \cite{GRALU1}, for algorithmic aspects see \cite{PAPRIN}
and ''non-classical'' applications can be found in \cite{PAGE1}, \cite{PAGE2}.
For a first  promising application of functional quantization to the pricing of financial derivatives through numerical integration on
path-spaces see
\cite{PAGE3}.

We address the issue of high-resolution quantization which concerns the performance of $n$-quantizers and the 
behaviour of $e_n(X)$ as $n \rightarrow \infty$.
The asymptotics of $e_n(X)$ for $\R^d$-valued random vectors has been completely elucidated for non-singular distributions $\P_X$ by the Zador 
Theorem (see \cite{GRALU1}) and for a class of self-similar (singular) distributions by \cite{GRALU2}. In infinite dimensions no such global
results hold, even for Gaussian processes.

It is convenient to use the symbols $\sim$ and $\lesim$, where $a_n \sim b_n$ means $a_n/b_n \rightarrow 1$ and $a_n \lesim b_n$ 
means $\limsup_{n \to \infty} a_n/b_n \leq 1$. A measurable function $\varphi : (s, \infty) \rightarrow (0,
\infty)\, (s \geq 0)$ is said to be regularly varying at infinity with index $b \in \R$ if, for every $c > 0$,
\[
\lim_{x \to \infty} \frac{\varphi(cx) }{\varphi(x)} = c^b .
\]

Now let $X$ be centered Gaussian.
Denote by $K_X \subset H$ the reproducing kernel Hilbert space (Cameron-Martin space) associated to the covariance operator
\begin{equation}
C_{_X} : H \rightarrow H, \; C_{_X} y := \E\,(< \!y, X\! >\!X)
\end{equation}
of $X$. Let $\lambda_1 \geq \lambda_2 \geq \ldots > 0$ be the ordered nonzero eigenvalues of $C_{_X}$ and let $\{u_j : j \geq 1 \}$ be 
the corresponding orthonormal basis of supp$(\P_X)$ consisting of eigenvectors (Karhunen-Lo\`eve basis). If
$d := \dim K_X < \infty$, then $e_n (X) = e_n \left( \bigotimes\limits^d_{j=1} N(0, \lambda_j)\right)$,
the minimal $n$th $L^2$-quantization error of
$\bigotimes\limits^d_{j=1} N(0, \lambda_j)$ with respect to the $l_2$-norm on $\R^d$, and thus we can read off the asymptotic
 behaviour of $e_n(X)$ from the high-resolution formula
\begin{equation}\label{1.7}
e_n (\bigotimes\limits^d_{j=1} N(0, \lambda_j)) \sim q(d) \sqrt{2 \pi} \left( \Pi^{d}_{j=1}
\lambda_j\right)^{1/2 d}
\left( \frac{d+2}{d}\right)^{(d+2)/4} n^{-1/d} \; \mbox{ as } \; n \rightarrow \infty
\end{equation}
where $q(d) \in (0, \infty)$ is a constant depending only on the dimension $d$ (see \cite{GRALU1}). Except in dimension $d=1$ and $d=2$, the
true value of $q(d)$ is unknown. However, one knows (see \cite{GRALU1}) that
\begin{equation}\label{(1.8)}
q(d) \sim \left( \frac{d}{2 \pi e}\right)^{1/2} \; \mbox{ as } \; d \rightarrow \infty.
\end{equation}

Assume $\dim K_X = \infty$. Under regular behaviour of the eigenvalues the sharp asymptotics of $e_n (X)$ can be derived
analogously to (\ref{1.7}). In view of~(\ref{(1.8)}) it is reasonable to expect that the limiting constants can be evaluated. The recent
high-resolution formula is as follows.

\begin{Thm}
(\cite{LUPA2}) Let $X$ be a centered Gaussian. Assume $\lambda_j \sim \varphi(j)$ as $j \rightarrow \infty$, where
$\varphi : (s, \infty) \rightarrow (0, \infty)$ is a decreasing, regularly varying function at infinity of 
index $-b < -1$ for some $s \geq 0$. Set, for every $x > s$,
\[
\psi (x) := \frac{1}{x \varphi(x)} .
\]
Then
\[
e_n(X) \sim \left(\left( \frac{b}{2}\right)^{b-1}\!\! \frac{b}{b-1}\right)^{1/2} \psi ( \log n )^{-1/2}
\;
\mbox{ as }
\;
\vspace{0.5cm} n \rightarrow \infty .
\]
\end{Thm}

A high-resolution formula in case $b = 1$ is also available (see \cite{LUPA2}). Note that the restriction $-b \leq -1$ 
on the index of $\varphi$ is natural since $\sum\limits^\infty_{j=1} 
\lambda_j < \infty$. The minimal $L^r$-quantization errors of $X$, $0 < r < \infty$, are strongly equivalent to the $L^2$-errors $e_n(X)$
(see \cite{DERE}) and thus exhibit the same high-resolution behaviour.

A related quantization problem is the Kolmogorov metric entropy problem for the closed unit ball
\begin{equation}
U_X := \left\{ x \in K_X : \|  x \|_{K_X} \leq 1 \} = \{ x \in \; \mbox{supp} (\P_X) :
\sum\limits_{j \geq 1} \frac{< \!x,u_j \!>^2}{\lambda_j} \leq 1 \right\}
\end{equation}
of $K_X$ (Strassen ball). Note that $U_X$ is a compact subset of $H$. For $n \in \N$, the metric entropy problem for $U_X$ consists in minimizing
\[
\max_{x \in U_X} \min_{a \in \alpha} \| x-a \| = \|\min_{a \in \alpha} \| X^{'} - a \|  \|_{L^\infty (
\P )}
\]
over all subsets $\alpha \subset H$ with $\mbox{card} (\alpha) \leq n$, where $X^{'}$ is any $H$-valued random
 vector with $\mbox{supp} (\P_{X^{'}} ) = U_X$.  The $n$th entropy number is then defined by
\begin{equation}
e_n(U_X) := \inf \left\{ \max_{x \in U_X} \min_{a \in \alpha} \| x-a \|  : \alpha \subset H, \;
\mbox{card}( \alpha ) \leq n \right\} .
\end{equation}

If $d := \dim K_X < \infty$, then $e_n(U_X) = e_n ({\cal E}_d)$, the $n$th entropy number of the ellipsoid
\[
{\cal E}_d := \left\{ x \in \R^d : \sum\limits^d_{j=1} \frac{x2_j}{\lambda_j} \leq 1 \right\}
\]
with respect to the $l_2$-norm on $\R^d$. Thus we can read off the asymptotic behaviour of $e_n (U_X)$
from the formula
\begin{equation}
e_n ({\cal E}) \sim p(d) ( \Pi^d_{j=1} \lambda_j)^{1/2} ( \mbox{vol} \; (B_d (0,1)))^{1/d} n^{-1/d} \; \mbox{ as } \; n \rightarrow \infty
\end{equation}
where the constant $p(d) \in (0,\infty)$ is unknown for $d \geq 3$ and $p(d) \sim q(d) , d \rightarrow \infty$ (see \cite{KOLMO}, \cite{GRALU1}).

If $\dim K_X = \infty$, the recent solution of the Kolmogorov metric entropy problem for $U_X$ is as follows.

\begin{Thm} (\cite{LUPA3})
Assume the situation of Theorem 1. Then
\[
e_n (U_X) \sim \left( \frac{b}{2}\right)^{b/2}\!\!\varphi(\log n)^{1/2} \; \mbox{ as } \; n \rightarrow
\infty .
\]
\end{Thm}

This formula is still valid for $b = 1$ and, ignoring the probabilistic interpretation, also for
$b \geq 0$ $(00 := 1)$ provided $\lambda_j \rightarrow 0$ as $j \rightarrow \infty$.
(see \cite{GRALU3}, \cite{LUPA3}). A different approach via the inverse of $e_n (U_X)$, the Kolmogorov
 $\varepsilon$-entropy, is due to Donoho \cite{DONO}. (However, his result does not provide the correct constant.)
>From Theorems~$1$ and~$2$ we conclude that functional quantization and metric entropy are related by
\begin{equation}
e_n (X) \sim \left( \frac{2 \log n}{b-1}\right)^{1/2}\!\!e_n (U_X) \; \mbox{ as } \; n \rightarrow \infty.
\end{equation}

The paper is organized as follows. In Section~2 we investigate Riemann-Liouville processes in
$H = L^2 ([0,1], dt)$. For $\rho \in (0, \infty)$, the
Riemann-Liouville process $X^\rho = (X^\rho_t)_{t \in [0,1]}$ on $[0,1]$ is defined by
\begin{equation}\label{(1.13)}
X^\rho_t := \int^t_0 (t - s)^{\rho- \frac{1}{2}} d W_s
\end{equation}
where $W$ is a standard Brownian motion.
We derive a high-resolution formula for $X^\rho$ and correspondingly, the precise entropy asymptotics for fractional Sobolev balls.
As a consequence we obtain a new result for fractionally integrated Brownian motions.
In Section 3 we describe a quantization procedure which furnishes asymptotically optimal quantizers in the situation of 
Theorem 1. Here the Karhunen-Lo\`eve expansion plays a crucial
r\^ole. In Section 4 we discuss a dimension conjecture.
\section{Riemann-Liouville processes}
\setcounter{equation}{0}
\setcounter{Assumption}{0}
\setcounter{Theorem}{0}
\setcounter{Proposition}{0}
\setcounter{Corollary}{0}
\setcounter{Lemma}{0}
\setcounter{Definition}{0}
\setcounter{Remark}{0}
Let $X^\rho = (X^\rho_t)_{t \in [0,1]}$ be the Riemann-Liouville process of index $\rho \in (0, \infty)$ as defined in~(\ref{(1.13)}). Its
covariance function is given by
\begin{equation}\label{(2.1)}
\E\,X^\rho_s X^\rho_t = \int^{s \wedge t}_0 (t - r)^{\rho- \frac{1}{2}} (s-r)^{\rho - \frac{1}{2}} dr .
\end{equation}
Using $\rho \wedge \frac{1}{2}$-H{\"o}lder continuity of the application
$t \mapsto X^\rho_t$ from [0,1] into $L^2(\P)$ and
the Kolmorogov criterion one checks that $X^\rho$ has a pathwise continuous modification so that we may assume without 
loss of generality that $X^\rho$ is pathwise continuous. In particular, $X^\rho$ can be seen as
a centered Gaussian random vector with values in
\[
H = L^2 ([0,1], dt).
\]
The following high-resolution formula relies on a theorem by Vu and Gorenflo \cite{TUAN} on singular values 
of Riemann-Liouville integral operators
\begin{equation}
R_\beta \,g(t) = \frac{1}{\Gamma(\beta)} \int^t_0 (t-s)^{\beta-1} g(s) ds , \hskip 1 cm \beta \in (0, \infty).
\end{equation}

\begin{Thm}
For every $\rho \in (0, \infty)$,
\[
e_n (X^\rho) \sim \pi^{-(\rho + \frac{1}{2})} (\rho + 1/2)^\rho ( \frac{2 \rho + 1 }{2 \rho} )^{1/2} \Gamma( \rho + 1/2) (\log n)^{-\rho} \;
\mbox{ as } \; n \rightarrow \infty .
\]
\end{Thm}

\noindent {\bf Proof.} For $\beta > 1/2$, the Riemann-Liouville fractional integral operator $R_\beta$
is a bounded operator from
$L^2 ([0,1], dt)$ into $ L^2 ([0,1], dt)$. The covariance operator
\[
C_\rho : L^2 ([0,1], dt) \rightarrow L^2 ([0,1], dt)
\]
of $X^\rho$ is given by the Fredholm transformation
\[
C_\rho g(t) = \int^1_0 g(s) E X^\rho_s X^\rho_t ds .
\]
Using~(\ref{(2.1)}), one checks that $C_\rho$ admits a factorization
\[
C_\rho = S_\rho S^{*}_\rho ,
\]
where
\[
S_\rho = \Gamma (\rho +1/2 ) R_{ \rho + \frac{1}{2}} .
\]
Consequently, it follows from Theorem 1 in \cite{TUAN} that the eigenvalues $\lambda_1 \geq \lambda_2 \geq \ldots > 0$ of $C_\rho$ satisfy
\begin{equation}\label{(2.3)}
\lambda_j \sim \Gamma ( \rho +1/2)^2 ( \pi j)^{-(2 \rho + 1)} \; \mbox{ as } \; j \rightarrow \infty.
\end{equation}
Now the assertion follows from Theorem 1 (with $\varphi(x) = \Gamma(\rho + 1/2)^2 \pi^{-b} x^{-b}$ and $b = 2 \rho + 1)$.
\hfill{$\Box$}

\bigskip
An immediate consequence for fractionally integrated Brownian motions on $[0,1]$ defined by
\begin{equation}
Y^\beta_t := \frac{1}{\Gamma(\beta)} \int^t_0 (t-s)^{\beta-1} W_s ds
\end{equation}
for $\beta \in (0, \infty)$ is as follows.

\begin{Cor}
For every $\beta \in (0, \infty),$
\[
e_n (Y^\beta) \sim \pi^{-(\beta+1)} (\beta+1)^{\beta + \frac{1}{2}}
 ( \frac{2 \beta + 2}{2 \beta +1})^{1/2} (\log n)^{-(\beta + \frac{1}{2})} \; \mbox{ as } \; n \to \infty .
\]
\end{Cor}

\noindent {\bf Proof.} For $\rho > 1/2$, the Ito formula yields
\[
X^\rho_t = (\rho - \frac{1}{2} ) \int^t_0 (t-s)^{\rho - \frac{3}{2}} W_s ds .
\]
Consequently,
\[
Y^\beta_t = \frac{1}{\beta \Gamma (\beta)} \beta \int^t_0 (t-s)^{\beta + \frac{1}{2} - \frac{3}{2} }
W_s ds = \frac{1}{\Gamma (1+\beta)} X_t^{\beta + \frac{1}{2} } .
\]
The assertion follows from Theorem 3. \hfill{$\Box$}

\bigskip
\noindent {\bf Remark.} The preceding corollary provides new high-resolution formulas for $e_n (Y^\beta)$ in the cases
$\beta \in (0, \infty) \setminus \N$.

\bigskip One  further consequence is a precise relationship between the quantization errors of Riemann-Liouville processes
and fractional Brownian motions. The fractional Brownian motion with Hurst exponent $\rho \in (0,1]$ 
is a centered pathwise continuous Gaussian process
$Z^\rho = (Z^\rho_t)_{t \in [0,1]}$ having the covariance function
\begin{equation}
\E\,Z^\rho_s Z^\rho_t = \frac{1}{2} ( s^{2 \rho} + t^{2 \rho} - \mid s-t \mid^{2 \rho} ) .
\end{equation}

\begin{Cor}
For every $\rho \in (0,1)$,
\[
e_n (X^\rho) \sim \frac{\Gamma(\rho + 1/2 )}{ ( \Gamma ( 2 \rho + 1) \sin ( \pi \rho))^{1/2} }
e_n (Z^\rho) \; \mbox{ as } \; n \rightarrow \infty .
\]
\end{Cor}
\bigskip
{\bf Proof.} By \cite{LUPA2}, we have
\[
e_n (Z^\rho) \sim \pi^{-(\rho + \frac{1}{2})} ( \rho + 1/2)^\rho \left( \frac{2 \rho +1}{2 \rho}\right)^{1/2} ( \Gamma (2
\rho +1)
\sin (\pi \rho))^{1/2} ( \log n)^{-\rho} , n \rightarrow \infty .
\]
Combining this formula with Theorem 3 yields the assertion \hfill{$\Box$}  

\bigskip
Observe that strong equivalence $e_n (X^\rho) \sim e_n (Z^\rho)$ as $n \rightarrow \infty$ 
is true for exactly two values of $\rho \in (0,1)$, namely for $\rho = 1/2$ where even $e_n (X^{1/2} ) = e_n (Z^{1/2} ) = e_n (W)$ and, a bit
mysterious, for $\rho = 0.81557\ldots$

Now consider the Strassen ball $U_\rho $ of $X^\rho$. Since the covariance operator $C_\rho$ satisfies
$C_\rho = \Gamma(\rho + \frac{1}{2}) R_{\rho + \frac{1}{2}} (\Gamma ( \rho + \frac{1}{2})
R_{\rho + \frac{1}{2}} )^{*}$,
one gets
\begin{eqnarray}
U_\rho & = & \Gamma ( \rho + 1/2) R_{\rho + \frac{1}{2}} (B_{L^2} (0,1)) \\
& = & \left\{ R_{\rho + 1/2} g: g \in L^2([0,1], dt), \int1_0 g(t)^2dt \leq \Gamma(\rho +
1/2)^2 \right\} , \nonumber
\end{eqnarray}
a fractional Sobolev ball. Theorem 2 and~(\ref{(2.3)}) yield the solution of the entropy problem for fractional Sobolev balls.
\begin{Thm}
For every $\rho \in (0, \infty)$,
\begin{eqnarray*}
e_n (U_\rho) & \sim & \left( \frac{\rho + \frac{1}{2} }{ \pi } \right)^{\rho + \frac{1}{2} }  \Gamma (\rho +1/2 )
         (\log n)^{-(\rho + \frac{1}{2})} \\
     & \sim & \left( \frac{\rho}{\log n}\right)^{1/2} e_n (X^\rho) \; \mbox{ as } \; n \rightarrow \infty .
\end{eqnarray*}
\end{Thm}
\section{Asymptotically optimal functional quantizers}
\setcounter{equation}{0}
\setcounter{Assumption}{0}
\setcounter{Theorem}{0}
\setcounter{Proposition}{0}
\setcounter{Corollary}{0}
\setcounter{Lemma}{0}
\setcounter{Definition}{0}
\setcounter{Remark}{0}
Let $X$ be a $H$-valued random vector satisfying~(\ref{(1.2)}). For every $n \in \N$, $L^2$-optimal $n$-quantizers $\alpha \subset H$ exist,
that is
\[
(\E\,\min_{a \in \alpha} \| X -a \|^2)^{1/2} = e_n (X) .
\]
If card (supp$(\P_X)) \geq n$, optimal $n$-quantizers $\alpha$ satisfy card$(\alpha) = n$,
$\P (X \in C_a (\alpha)) > 0$ and the stationarity condition
\[
a = \E\,(X \mid  \{ X \in C_a (\alpha)\}), \,a \in \alpha
\]
or what is the same
\begin{equation}\label{(3.1)}
\hat{X}^\alpha = \E\,(X \mid  \hat{X}^\alpha)
\end{equation}
for every Voronoi partition $\{ C_a(\alpha) : a \in \alpha \}$ (see \cite{LUPA1}). In particular,
$\E\,\hat{X}^\alpha = \E\,X$.

Now let $X$ be centered Gaussian with $\dim K_X = \infty$. The Karhunen-Lo\`eve basis $\{u_j : j \geq 1 \}$
consisting of normalized eigenvectors of $C_{_X}$ is optimal for the quantization of Gaussian random vectors
(see \cite{LUPA1}). So we start with the Karhunen-Lo\`eve expansion
\[
X \stackrel{H}{=} \sum\limits^\infty_{j=1} \lambda^{1/2}_j Z_j u_j ,
\]
where $Z_j = < \!X, u_j\!> / \lambda^{1/2}_j, j \geq 1$ are i.i.d. $N(0,1)$-distributed random variables. The design of an
asymptotically optimal quantization of $X$ is based on optimal quantizing blocks of coefficients of variable ($n$-dependent)
block length. Let $n \in \N$ and fix temporarily
$m, l, n_1, \ldots , n_m \in \N$ with $\Pi^m_{j=1} n_j \leq n$, where $m$ denotes the number of blocks, $l$ 
the block length and $n_j$ the size of the quantizer for the $j$th block
\[
Z^{(j)} := (Z_{(j-1)l+1} , \ldots , Z_{jl}), \quad j \in \{ 1, \ldots , m\}.
\]
Let $\alpha_j \subset \R^l$ be an $L^2$-optimal $n_j$-quantizer for
$Z^{(j)}$ and let $\widehat{Z^{(j)}} =  \widehat{Z^{(j)}}^{\alpha_j}$
be a $\alpha_j$-quantization of $Z^{(j)}$. Then, define a quantized version of $X$ by
\begin{equation}\label{3.2}
\hat{X}^n := \sum\limits^m_{j=1} \sum\limits^l_{k=1} \lambda^{1/2}_{(j-1)l+k} ( \widehat{Z^{(j)}} )_k
u_{(j-1)l+k} .
\end{equation}
It is clear that card$(\hat{X}^n(\Omega )) \leq n$. Using (\ref{(3.1)}) for $Z^{(j)}$, one gets
$\E\,\hat{X}^n = 0$. If
\[
\widehat{Z^{(j)}} = \sum\limits_{b \in \alpha_j} b \mbox{\bf 1}_{C_b(\alpha_j)} (Z^{(j)} ),
\]
then
\[
\hat{X}^n = \sum\limits_{a \in \times^m_{j=1} \alpha_j} ( \sum\limits^m_{j=1} \sum\limits^l_{k=1}
\lambda^{1/2}_{(j-1) l+k} a^{(j)}_k u_{(j-1) l+k} ) \Pi^m_{j=1} \mbox{\bf 1}_{C_{a^{(j)}} (\alpha_j)} (Z^{(j)} )
\]
where $a = (a^{(1)}, \ldots , a^{(m)} ) \in \times^{m}_{j=1} \alpha_j$. Observe that in general,
$\hat{X}^n$ is not a Voronoi quantization of $X$ since it is based on the (less complicated) Voronoi 
partitions for $Z^{(j)}, j \leq m$. $(\hat{X}^n$ is a Voronoi quantization if $l = 1$ or if
$\lambda_{(j-1) l+1} = \ldots = \lambda_{jl}$ for every $j$.)
Using again (\ref{(3.1)}) for $Z^{(j)}$ and the independence structure, one checks that $\hat{X}^n$
satisfies a kind of stationarity equation:
\[
\E\,(X \mid  \hat{X}^n) = \hat{X}^n .
\]
\begin{Lem} Let $n\ge1$. For every $l\ge 1$ and every $m\ge 1$
\begin{equation}\label{Pythagore}
\E\,\| X - \hat{X}^n \|^2 \leq \sum\limits^m_{j=1} \lambda_{(j-1)l+1} e_{n_{j}} (N(0, I_l))^2
+ \sum\limits_{j \geq ml+1} \lambda_j .
\end{equation}
Furthermore,~(\ref{Pythagore}) stands as an equality if $l=1$ (or $\lambda_{(j-1) l+1} = \ldots = \lambda_{jl}$ for every $j,\,l\ge 1$).
\end{Lem}

\noindent{\bf Proof.}
The claim follows from the orthonormality of the basis $\{ u_j : j \geq 1 \}$. We have
\[
\begin{array}{lcl}
\E\,\| X - \hat{X}^n \|^2
& = & \sum\limits^m_{j=1} \sum\limits^l_{k=1} \lambda_{(j-1) l+k} \E\,\mid Z^{(j)}_k - (\widehat{Z^{(j)}})_k
        \mid^2 + \sum\limits_{j \geq ml +1} \lambda_j \\
& \leq & \sum\limits^m_{j=1} \lambda_{(j-1)l+1} \sum\limits^l_{k=1} \E\,\mid Z^{(j)}_k - \widehat{Z^{(j)}} )_k \mid^2 +
        \sum\limits_{j \geq ml +1} \lambda_j \\
& = & \sum\limits^m_{j=1} \lambda_{(j-1)l+1} e_{n_{j}} (Z^{(j)})^2 + \sum\limits_{j \geq ml +1} \lambda_j .
\end{array}
\]
\hspace*{\fill}{$\Box$}

\bigskip
Set
\begin{equation}\label{3.3}
C(l) := \sup_{k \geq 1} k^{2/l} e_k (N(0,I_l))^2.
\end{equation}
By~(\ref{1.7}), $C(l) < \infty$. For every $l \in \N$,
\begin{equation}\label{(3.4)}
e_{n_{j}} (N(0,I_l)^2 \leq n^{-2/l}_j C(l)
\end{equation}
Then one may replace the optimization problem which consists, for fixed $n$, in 
minimizing the right hand side of Lemma 1 by the following optimal allocation problem:
\begin{equation}\label{optipb}
\min\{ C(l) \sum\limits^m_{j=1} \lambda_{(j-1)l+1} n^{-2/l}_j + 
\sum\limits_{j \geq ml +1} \lambda_j : m, l, n_1 , \ldots , n_m \in \N, \Pi^m_{j=1} n_j \leq n \} .
\end{equation}
Set
\begin{equation}\label{mln}
m = m(n,l) := \max \{ k \geq 1 : n^{1/k} \lambda_{(k-1)l+1}^{l/2} ( \Pi^k_{j=1} \lambda_{(j-1)l+1} )^{-l/2k}
\geq 1 \} ,
\end{equation}\label{nj}
\begin{equation}
n_j = n_j (n,l) := [n^{1/m} \lambda^{l/2}_{(j-1)l+1} ( \Pi^m_{i=1} \lambda_{(i -1)l+1} )^{-l/2m}] ,\; j \in \{ 1 , \ldots , m\},
\end{equation}
where $[x]$ denotes the integer part of $x \in \R$ and
\begin{equation}\label{ln}
l = l_n := [(\max \{ 1 , \log n \})^\vartheta ], \;\vartheta \in (0,1) .
\end{equation}
In the following theorem it is demonstrated that this
choice is at least asymptotically optimal provided the eigenvalues are regularly varying.

\begin{Thm}
Assume the situation of Theorem 1. Consider $\hat{X}^n$ with tuning parameters defined in~(\ref{mln})-(\ref{ln}).
Then $\hat{X}^n$ is asymptotically $n$-optimal, i.e.
\end{Thm}
\[
(\E\,\| X - \hat{X}^n \|^2 )^{1/2} \sim e_n (X) \; \mbox{ as } \; n \rightarrow \infty .
\]

Note that no block quantizer with fixed block length is asymptotically optimal (see \cite{LUPA2}). As mentioned above,
$\hat{X}^n$ is not a Voronoi quantization of $X$. If $\alpha_n := \hat{X}^n(\Omega)$, then the Voronoi
quantization $\hat{X}^{\alpha_n}$ is clearly also asymptotically $n$-optimal.

\bigskip
The key property for the proof is the following $l$-asymptotics of the constants $C(l)$ defined in~(\ref{3.3}). It is
interesting to consider also the smaller constants
\begin{equation}
Q(l) := \lim_{k \to \infty} k^{2/l} e_k (N(0, I_l))^2
\end{equation}
(see (\ref{1.7})).

\begin{Pro}
The sequences $(C(l))_{l \geq 1}$ and $(Q(l))_{l \geq 1}$ satisfy
\[
\lim_{l \to \infty} \frac{C(l)}{l} = \lim_{l \to \infty} \frac{Q(l)}{l} = \inf_{l \geq 1} \frac{C(l)}{l}
= \inf_{l \geq 1} \frac{Q(l)}{l} = 1.
\]
\end{Pro}
{\bf Proof.} From \cite{LUPA2} it is known that
\begin{equation}
\liminf_{l \to \infty} \frac{C(l)}{l} = 1.
\end{equation}
Furthermore, it follows immediately from~(\ref{1.7}) and~(\ref{(1.8)}) that
\begin{equation}
\lim_{l \to \infty} \frac{Q(l)}{l} = 1.
\end{equation}
(The proof of the existence of $\displaystyle \lim_{l \to \infty} C(l)/l$ we owe to S. Dereich.) For $l_0, l \in \N$ with
$l \geq l_0$, write
\[
l = n\, l_0 + m \; \mbox{with} \; n \in \N, m \in \{ 0, \ldots , l_0 -1 \} .
\]
Since for every $k \in \N$,
\[
[k^{l_0/l}]^n \; [k^{1/l}]^m \leq k ,
\]
one obtains by a block-quantizer design consisting of $n$ blocks of length $l_0$ and $m$ blocks of length 1 for
quantizing $N(0, I_l)$,
\begin{equation}\label{(3.13)}
e_k (N(0,I_l))^2 \leq n e_{[k^{l_0/l}]} (N(0,I_{l_0}))^2 + m e_{[k^{1/l}]} (N(0,1))^2 .
\end{equation}
This implies

\begin{eqnarray*}
C(l) &  \leq & n C(l_0)  \sup_{k \geq 1} \frac{k^{2/l}}{[k^{l_0/l}]^{2/l_0}} + m C(1) \sup_{k \geq 1} \frac{k^{2/l}}{[k^{1/l}]^2}  \\
     & \leq    & 4^{1/l_0} n C(l_0) + 4 m C(1) .
\end{eqnarray*}

Consequently, using $n/l \leq 1/l_0$,
\[
\frac{C(l)}{l} \leq \frac{ 4^{1/l_0} C(l_0) }{ l_0 } + \frac{ 4 m C(1) }{ l }
\]
and hence
\[
\limsup_{l \to \infty} \frac{C(l)}{l} \leq \frac{ 4^{1/l_0} C(l_0) }{ l_0 } .
\]
This yields
\begin{equation}
\limsup_{l \to \infty} \frac{C(l)}{l} \leq \liminf_{l_0 \rightarrow \infty} \frac{C(l_0)}{l_0} = 1 .
\end{equation}
It follows from~(\ref{(3.13)}) that
\[
Q(l) \leq n Q(l_0) + m Q(1) .
\]
Consequently
\[
\frac{Q(l)}{l} \leq \frac{Q(l_0)}{l_0} + \frac{m Q(1)}{l}
\]
and therefore
\[
1 = \lim_{l \to \infty} \frac{Q(l)}{l} \leq \frac{ Q(l_0) }{ l_0 } .
\]
This implies
\begin{equation}
\inf_{l_0 \geq 1} \frac{Q(l_0)}{l_0} = 1.
\end{equation}
Since $Q(l) \leq C(l)$, the proof is complete. \hfill{$\Box$}

\bigskip
The $n$-asymptotics of  the number $m(n, l_n)l_n$ of quantized coefficients in the Karhunen-Lo\`eve expansion  in
the quantization $\hat{X}^n$ is as follows.
\begin{Lem}
(\cite{LUPA3}, Lemma 4.8) Assume the situation of Theorem 1. Let $m(n, l_n)$ be defined by~(\ref{mln}) and~(\ref{ln}). Then
\[
m(n, l_n) l_n \sim \frac{2 \log n}{b} \; \mbox{ as } \; n \rightarrow \infty .
\]
\end{Lem}

\noindent {\bf Proof of Theorem 5.} For every $n \in \N$,
\begin{eqnarray*}
\sum\limits^m_{j=1} \lambda_{(j-1)l+1} n^{-2/l}_j
& \leq & \sum\limits^m_{j=1} \lambda_{(j-1)l+1} (n_j +1)^{-2/l} ( \frac{n_j +1}{n_j} )^{2/l} \\
& \leq & 4^{1/l} m n^{-2/ml} ( \Pi^m_{j=1} \lambda_{(j-1)l+1} )^{1/m} \\
& \leq & 4^{1/l} m \lambda_{(m-1) l+1} .
\end{eqnarray*}
Therefore, by Lemma 1 and~(\ref{(3.4)}),
\[
\E\,\| X - \hat{X}^n \|^2 \leq 4^{1/l} \frac{C(l)}{l} m l \lambda_{(m-1)l+1} + \sum\limits_{j \geq m l + 1} \lambda_j
\]
for every $n \in \N$. By Lemma~2, we have
\[
m l = m(n,l_n) l_n \sim \frac{2 \log n}{b} \; \mbox{ as } \; n \rightarrow \infty .
\]
Consequently, using regular variation at infinity with index $-b < -1$ of the function $\varphi$,
\[
m l \lambda_{(m-1)l+1} \sim m l \lambda_{ml} \sim \left( \frac{2}{b} \right)^{1-b} \psi (\log n)^{-1}
\]
and
\[
\sum\limits_{j \geq m l +1} \lambda_j \sim \frac{m l \varphi(ml)}{b-1} \sim
\frac{1}{b-1} \left( \frac{2}{b} \right)^{1-b} \psi ( \log n)^{-1} \; \mbox{ as } \; n \rightarrow \infty,
\]
where, like in Theorem 1, $\psi(x) = 1/x \varphi(x)$. Since by Proposition 1,
\[
\lim_{n \to \infty} \frac{4^{1/l_n}C(l_n)}{l_n} = 1 ,
\]
one concludes
\[
\E\,\| X - \hat{X}^n \|^2 \stackrel{<}{\sim} \left( \frac{2}{b} \right)^{1-b} \frac{b}{b-1}\, \psi ( \log n)^{-1} \; \mbox{ as }
\; n \rightarrow \infty .
\]
The assertion follows from Theorem 1. \hfill{$\Box$}

\bigskip
 \noindent {\sc Numerical and computational aspects: } As soon as the Karhunen-Lo\`eve  basis $(u_j)_{j\ge 1}$ of a Gaussian
process
$X$ is explicit, it is possible to compute the asymptotically optimal 
functional   quantization~(\ref{3.2})  which solves the
minimization problem~(\ref{optipb}) as well as its distribution and 
induced quantization error (at least for a given $\vartheta\!\in(0,1)$).
This is   possible since  some  optimal (or at least locally optimal) 
vector quantizations of the $N(0, I_d)$-distribution has been 
already computed and kept off line. Let  us be more specific.

\smallskip
-- In
$1$-dimension, the normal distribution $N(0, 1)$ has only one
stationary
$n$-quantizer --  hence optimal -- since its probability density is $\log$-concave (for this result due to Kiefer, see
$e.g.$~\cite{GRALU1}). Deterministic methods to compute these optimal quantizers are based on the
stationary equation~(\ref{(3.1)}). They are very easy to implement, converge very
fast with a very high accuracy. The Newton-Raphson algorithm is a
possible choice (see~\cite{PAPRIN} for details). Closed forms for the lowest
quadratic quantization error $\|Z-\widehat Z\|_{L^2(\P)}$ and  for the
distribution of the optimal $n$-quantization $\widehat Z^\alpha$ as a function of the optimal
$n$-quantizer $\alpha$      are also available in~\cite{PAPRIN}. 
These three quantities have been tabulated up to very high values of $n$. A file can be downloaded at
the URL {\tt  www.proba.jussieu.fr/pageperso/pages.html}.

\smallskip
-- In higher dimension, one still relies on the  stationary
equation (3.1) which reads: 
\[
\E\left(\mbox{\bf 1}_{C_a(\alpha)}(Z) (a-Z)\right)= 0,\qquad a\!\in\alpha
\]
 One must keep in mind that the left hand side of the above equation is
but the gradient of the (squared) quantization error 
$\E \| Z- \widehat Z ^\alpha \|^2$ viewed as a function of the
quantizer $\alpha$ (assumed to be of full size $n$). A stochastic
gradient  descent based on this integral representation can be
implemented easily since the normal distribution $N (0,I_d)$ can
be  simulated on a computer from (pseudo-)random numbers ($e.g.$ by the
Box-Muller method). This algorithm is known as  the Competitive Learning
Vector Quantization (or $CLVQ$) algorithm. It has been extensively
investigated   both from a theoretical   (see $e.g.$~\cite{PAGE1}, \cite{BOPA}) and numerical 
(see $e.g.$~\cite{PAPRIN} as concerns normally distributed vectors) viewpoints. The algorithm reads
as follows: let
$(\zeta(t))_{t\ge 1}$ be an i.i.d. sequence of $N (0, I_d)$-distributed
random vectors, let $(\gamma_t)_{t\ge 1}$ be a decreasing sequence of
positive {\em gain} parameter satisfying $\sum_t \gamma_t=+\infty$ and $\sum_{t\ge
1} \gamma^2_t<+\infty$ and let
$\alpha(0)\!\in (\R^d)^n$ denote a starting 
$n$-quantizer. Then, at time $t\!\in\N$, one update the running
$n$-quantizer $\alpha(t-1):=(\alpha_1(t-1),\ldots,\alpha_n(t-1))$ as
follows
\begin{eqnarray*}
\mbox{\sc Competitive phase:} \;\; {\rm select } \quad i(t) &\in& {\rm argmin}
\{i : \; \|\alpha_i(t-1)-\zeta(t)\|= \min_j\|\alpha_j(t-1)-\zeta(t)\|\}\;\;\;\;
\\
 \mbox{\sc Learning phase: }\qquad\; \alpha_{i(t)} (t-1)& = &(1-\g_t)
\alpha_{i(t)}(t-1)+\g_t\,\zeta(t)\;\;\\
\hskip 1,5 cm \alpha(t)_j&=& \alpha_{j-1}(t-1),\qquad j\neq i(t).\;\;
\end{eqnarray*}
Some further details concerning the numerical implementation of this
  procedure can be found in \cite{PAPRIN}, especially some heuristics concerning
the initialization and the specification of the gain parameter sequence usually choosen of the form
$\gamma_t = \frac{A}{B+t}$. It converges toward some local minima of the
quantization error at a $\sqrt{\gamma_t}$-rate. Some $d$-dimensional grids ($d=2$ up to $10$) can be downloaded at the above URL for many
values of $n$ in the range $2$ up to $2\,000$. These quantizations were carried out to
solve numerically multi-dimensional stopping time problems (pricing of
American options on baskets, see \cite{PAGE2} and the references therein).

\medskip
The $1$-dimensional optimal quantization of the $N (0,1)$-distribution has already been used to produce some  optimal {\em scalar}
product functional quantization - $i.e.$ based on blocks of fixed length 1-
in~\cite{PAGE3} with some promising applications
to the pricing of path-dependent European options in stochastic volatility models (this work is also  based on  results about
 diffusion processes from~\cite{LUPA4}).   To be competitive with other methods (Monte Carlo, pde's) one needs to have good performances
for not too large values of $n$. Within this range of values, it is more efficient to perform directly a numerical optimisation of (3.3)
(or~(\ref{optipb})) with $l=1$ rather than using the theoretical asymptotically optimal parameters (3.7) and~(3.8).

As far as numerical implementation of functional quantization with $n$-varying block
size is concerned, some first numerical experiments carried out by Benedikt Wilbertz \cite{WILB} 
suggest that it slightly improves the scalar approach for high values of
$n$, say $n\le 10^6$, simply using up to $3$-dimensional $n_j$-quantizers  with some
$n_j$ not greater than $100$. A similar improvement can vbe obtained for lower values of $n$ (say $n\ge 20\,000$) by using product quantizers
made of blocks with mixed diemnsions ($1$, $2$ or $3$).

\bigskip
\noindent {\sc Examples:} The basic example (among Riemann-Liouville processes) is $X^{1/2} = W$ and $H = L^2 ([0,1], dt)$, where
\begin{equation}
\lambda_j = ( \pi(j- 1/2))^{-2} , \;u_j(t) = \sqrt{2} \; \mbox{sin} \, \left(t/\sqrt{\lambda_j} \right),\;  j \geq 1.
\end{equation}
Since for $\delta, \rho \in (0, \infty)$,
\[
X^{\delta + \rho} = \frac{\Gamma(\delta + \rho + \frac{1}{2})}{\Gamma( \rho + \frac{1}{2})} R_\delta (X^\rho) ,
\]
one gets expansions of $X^{\delta + \rho}$ from Karhunen-Lo\`eve expansions of $X^\rho$. In particular,
\[
X^{\delta + \frac{1}{2}} = \Gamma (\delta+1) \sum\limits^\infty_{j=1} \sqrt{\lambda_j} Z_j R_\delta (u_j).
\]
However, the functions $R_\delta(u_j), j \geq 1$, are not orthogonal in $H$ so that the nonzero 
correlation between the components of $(Z^{(j)} - \widehat{Z^{(j)}})$ prevents the previous estimates for $\E \|X-\widehat X^n\|^2$  given in
Lemma~1 from working in this setting in the general case.

However, when $l=1$ (scalar product quantizers made up with blocks of fixed length
$l=1$), one checks that these estimates still stand as equalities since orthogonality can now be substituted by 
the independence of $Z_j - \hat{Z}_j$ 
and stationarity property~(\ref{(3.1)}) of the quantizations $\hat{Z}_j, j \geq 1$.
It is often good enough for applications
to use  scalar product quantizers (see
\cite{LUPA1},
\cite{PAGE3}). If, for instance $\delta = 1$, then
\[
X := X^{3/2} = \sum\limits^\infty_{j=1} \sqrt{\lambda_j} Z_j R_1 (u_j)  ,
\]
where
\[
R_1 (u_j) (t) =\sqrt{2 \lambda_j} (1 - \cos(t / \sqrt{\lambda_j})).
\]
Note that  $\displaystyle \| R_1 (u_j) \|^2 = \lambda_j (3 - 4 (-1)^{j-1}\sqrt{\lambda_j} ),\; j \geq 1$.
Set
\[
\hat{X}^n = \sum\limits^m_{j=1} \sqrt{\lambda_j} \hat{Z}_j R_1 (u_j) .
\]
 The quantization $\widehat X^n$ is non Voronoi  (it is
related to the Voronoi tessellation  of
$W$) and satisfies
\begin{equation}
\E\|X-\widehat X^n\|^2= \sum^m_{j = 1}  \lambda2_j (3 - 4 (-1)^{j-1}\sqrt{\lambda_j} )
e_{n_j} (N(0, 1))^2 + \sum_{j \geq m +1} \lambda^2_j (3-4(-1)^{j-1} \sqrt{\lambda_j}).
\end{equation}
It is possible to optimize the (scalar product) quantization error using this expression instead of~(\ref{optipb}).
As concerns asymptotics, if the  parameters  are tuned following~(\ref{mln})-(\ref{ln}) with $l=1$ and
$\lambda_j$ replaced by
\[
\nu_j := \lambda^2_j (3 + 4 \sqrt{\lambda_j} ) \sim 3 \pi^{-4} j^{-4} \quad\mbox{ as }\quad n\to \infty,
\]
and using Theorem 3 gives
\begin{equation}
( \E\,\| X - \hat{X}^n \|^2)^{1/2} \stackrel{<}{\sim} \left( \frac{3(12 C(1)+1)}{4}\right)^{1/2} e_n (X) \; \mbox{ as } \; n \rightarrow
\infty .
\end{equation}
Numerical experiments seem to confirm that $C(1) = Q(1)$. Since $Q(1) = \pi \sqrt{3}/2$ (see \cite{GRALU1}, p. 124), the 
above upper bound is then
\[
\left( \frac{3(6 \pi \sqrt{3} +1)}{4}\right)^{1/2} = 5.02357 \ldots
\]
\section{Dimension}
\setcounter{equation}{0}
\setcounter{Assumption}{0}
\setcounter{Theorem}{0}
\setcounter{Proposition}{0}
\setcounter{Corollary}{0}
\setcounter{Lemma}{0}
\setcounter{Definition}{0}
\setcounter{Remark}{0}
Let $X$ be a $H$-valued random vector satisfying~(\ref{(1.2)}). For $n \in \N$, let ${\cal C}_n(X)$ be the (nonempty) 
set of all $L^2$-optimal
$n$-quantizers. Introduce the integral number
\begin{equation}
d_n(X) := \min \left\{ \mbox{dim} \; \mbox{span} \; (\alpha) ) : \alpha \in {\cal C}_n (X) \right\} .
\end{equation}
It represents the dimension at level $n$ of the functional quantization problem for $X$. Here span$(\alpha)$ denotes the linear
subspace of $H$ spanned by $\alpha$. In view of Section 3, a reasonable conjecture for Gaussian random vectors is $d_n(X) \sim 2 \log n/b$ in regular
cases, where $-b$ is the regularity index. We have at least the following lower estimate in the Gaussian case.
\begin{Pro}
Assume the situation of Theorem 1. Then
\[
d_n (X) \stackrel{>}{\sim} \frac{1}{b^{1/(b-1)}} \; \frac{2 \log n}{b} \; \mbox{ as } \;
n \rightarrow \infty . \vspace{0.5cm}
\]
\end{Pro}
{\bf Proof.} For every $n \in \N$, we have
\begin{equation}
d_n (X) = \min \left\{ k \geq 0 : e_n ( \bigotimes^k_{j=1} N(0, \lambda_j))^2 + \sum\limits_{j \geq k+1} \lambda_j \leq e_n (X)^2 \right\}
\end{equation}
(see \cite{LUPA1}). Define
\[
c_n := \min \left\{ k \geq 0 : \sum\limits_{j \geq k + 1} \lambda_j \leq e_n (X)^2 \right\} .
\]
Clearly, $c_n$ increases to infinity as $n \rightarrow \infty$ and by (4.2), $c_n \leq d_n (X)$ for every $n \in \N$. 
Using Theorem 1 and the fact that $\psi$ is regularly varying at infinity with index $b-1$, we obtain
\[
((b-1) \psi(c_n))^{-1}  \sim  \sum\limits_{j \geq c_n +1} \lambda_j \sim e2_n(X)
 \sim  \left( \frac{2}{b} \right)^{1-b} \frac{b}{b-1} \,\psi (\log n)^{-1}
\]
and thus
\[
\psi(c_n)  \sim \left( \frac{2}{b} \right)^{1-b} \frac{1}{b} \psi (\log n)
\sim  \psi  \left( \frac{1}{b^{1/(b-1)}} \; \frac{2 \log n}{b}\right) \; \mbox{ as } \; n \rightarrow \infty .
\]
Consequently,
\[
c_n \sim \frac{1}{b^{1/(b-1)}} \; \frac{2 \log n}{b} \; \mbox{ as } \; n \rightarrow \infty .
\]
This yields the assertion. \hfill{$\Box$}  

\bigskip
For Riemann-Liouville processes one concludes
\[
d_n(X^\rho) \stackrel{>}{\sim} (2 \rho + 1)^{-1/2 \rho}  \; \frac{2 \log n}{2 \rho + 1}
\]
(see (2.3)).  

\bigskip For the metric entropy problem one may introduce the numbers $d_n (U_X)$ analogously. Then,
in the situation of Theorem 1 it is known that $d_n (U_X) \stackrel{>}{\sim} 2 \log n/b$ (see \cite{LUPA3}). 
It remains an open question whether $d_n (X)
\sim d_n(U_X) \sim 2 \log n/b$.

\end{document}